\documentclass[12pt,a4paper]{article}
\usepackage[dvipdfmx,hiresbb]{graphicx}
\usepackage{latexsym}
\usepackage{amssymb}
\usepackage{amsmath}
\usepackage{amscd}
\usepackage{mathrsfs}
\usepackage{bm}

\newtheorem{Th}{Theorem}[section]
\newtheorem{Co}[Th]{Corollary}
\newtheorem{Lem}[Th]{Lemma}

\newtheorem{Rem}[Th]{Remark}

\newtheorem{Pro}[Th]{Proposition}
\newcommand{\demo}{\par\noindent{\it Proof. \/}\ }
\newcommand{\enD}{\hfill $\Box$ \vspace{3truemm}\par}
\newcommand{\bx}{\mbox{\boldmath $x$}}
\newcommand{\bX}{\mbox{\boldmath $X$}}
\newcommand{\be}{\mbox{\boldmath $e$}}

\newcommand{\bv}{\mbox{\boldmath $v$}}
\newcommand{\by}{\mbox{\boldmath $y$}}
\newcommand{\bo}{\mbox{\boldmath $0$}}

\newcommand{\bn}{\mbox{\boldmath $n$}}

\newcommand{\blambda}{\mbox{\boldmath $\lambda$}}

\newcommand{\R}{{\mathbb R}}
\newcommand{\lon}{\longrightarrow}

\newcommand{\ou}{{\overline{u}}}

\begin{document}

\title{ Caustics of world hyper-sheets in \\ the Minkowski space-time}

\author{Shyuichi IZUMIYA}

\date{\today}

\maketitle

\begin{abstract}
 In the Minkowski space-time, 
a world hyper-sheet is a timelike hypersurface consisting of a one-parameter family of spacelike submanifolds.
Recently, Bousso and Randall introduced the notion of caustics of world hyper-sheets in order to define the notion of holographic domains in space-time.
Here, we give a mathematical framework for describing the caustics of world hyper-sheets in the Minkowski space-time.
As a consequence, we investigate the singularities of the caustics of world hyper-sheets and whose geometrical meanings.
Although the Minkowski space-time has zero gravity, this framework gives a simple toy model for general cases.
\end{abstract}
\renewcommand{\thefootnote}{\fnsymbol{footnote}}
\footnote[0]{2010 Mathematics Subject classification. Primary 58K05,57R45,32S05 ; Secondary 58K25, 58K60}
\footnote[0]{Keywords. Wave front propagations, Caustics}
\section{Introduction}
In this paper we investigate geometrical properties of caustics of world hyper-sheets in the Minkowski space-time as an application of the theory of graph-like Legendrian unfoldings \cite{Graph-like}.
Caustics appear in several area in Physics (i.e. geometrical optics \cite{Nye}, the theory of underwater acoustics \cite{Brekhov} and the theory of gravitational lensings\cite{Peters}
, and so on) and Mathematics (i.e. classical differential geometry \cite{Porteous} and theory of differential equations \cite{Hormander,IzuHJ93}, and so on \cite{Arnold-pink}).
The notion of caustics originally belongs to geometrical optics. We can observe the caustic formed by the rays reflected at a mirror. One of the examples of caustics in the classical differential geometry is the evolute of a curve in the Euclidean 
plane which is given by the envelope of normal lines emanated from the curve.
The ray in the Euclidean plane is considered to be a line, so that the evolute is the caustic in the sense of geometrical optics.  Moreover, the singular points of the evolute correspond to the vertices of the original curve. The vertex is the point at where the curve has higher order contact with the osculating circle (i.e. the point where the curvature has an extremum).
Therefore, the evolute provides an important geometrical information of the curve. We have the notion of evolutes for general hypersurfaces in the Euclidean space similar to
the plane curve case. In particular, there are detailed investigations on evolutes for surfaces in Euclidean $3$-space \cite {Izu07,Porteous}.
\par
The evolute of a hypersurface can be defined in the Minkowski space-time analogous to the Euclidean case. If we consider a timelike hypersurface in the Minkowski space-time, the normal line is
directed by a spacelike vector, whose speed exceeds the speed of the ray.
Therefore the evolute of a timelike hypersurface is not a caustic in the sense of Physics.
In the Minkowski space-time, the ray emanate from a codimension two spacelike submanifold is a normal line of the submanifold whose directer vector is lightlike, so the family of rays forms a lightlike hypersurface (i.e. a {\it light-sheet}).
The set of critical values of the light sheet is called a {\it lightlike focal set} along the spacelike submanifold.
Actually, the notion of light-sheets plays an important role in Physics which provides models of several kinds of horizons in space-times \cite{Ch}.
\par
On the other hand, a {\it world hyper-sheet} in the Minkowski space-time is a timelike hypersurface consisting of a one-parameter family of spacelike submanifolds of codimension two in the ambient space. Each spacelike submanifold is called a {\it momentary space}.
We consider the family of lightlike hypersurfaces along monetary spaces in the world hyper-sheet.
In \cite{Bousso, Bousso-Randall}, Bousso and Randall considered that
the locus of the singularities (the lightlike focal sets) of lightlike hypersurfaces along momentary spaces form a caustic in the Minkowski space-time.
This construction is originally from the theoretical physics (the string theory, the brane world scenario, the cosmology, and so on).
We call it a {\it BR-caustic} of the world sheet.
Moreover, we have no notion of the time constant in the relativity theory. Hence everything that is moving depends on the time.
Therefore, we consider world hyper-sheets in the relativity theory.
\par
In this paper we investigate the geometrical properties of BR-caustics as an application of the theory of graph-like Legendrian unfoldings \cite{Graph-like}.
In \S 2 we give the basic notions related to the Minkowski space-time.
Basic geometrical frame work for world hyper-sheets is given in \S 3.
The light-sheet along a momentary space is introduced in \S 4 and some calculations are given by using Lorentz distance squared functions.
In \S 5 the calculations in \S 4 are interpreted from the view point of contact with lightcones.
We briefly review the theory of graph-like Legendrian unfoldings in \S 6.
The notion of unfolded lightcone focal sets is introduced as a special case of the graph-like wave front in \S 7.
The caustic and the Maxwell set of the graph-like wave front are naturally induced.
In \S 8 the BR-caustic and the BR-Maxwell set are defined as the caustic and the Maxwell set of the graph-like wave front with respect to the distance squared function.
We give a classification of the caustics of world sheets in the $3$-dimensional Minkowski space-time in \S 9.
As a consequence, the local classification of BR-cautics in \S 9 is different from the local classification of the evolutes of timelike surfaces in the $3$-dimensional Minkowski space-time.

%%%%%%%%%%%%%%%%%%%%%%%%%%%%%%%%%
\section{The Minkowski space-time}
%%%%%%%%%%%%%%%%%%%%%%%%%%%%%%%%%%
\par 
We now introduce some basic notions on the $(n+1)$-dimensional Minkowski
space-time. For basic
concepts and properties, see \cite{Oneil}.
Let $\R^{n+1}=\{(x_0,x_1,\dots  ,x_n)\ |\ x_i \in \R\ (i=0,1,\dots , n)\ \}$
be an $(n+1)$-dimensional cartesian space. For any
$\bx =(x_0,x_1,\dots ,x_n),\ \by =(y_0,y_1,\dots ,y_n)\in \R^{n+1},$
the {\it pseudo scalar product } of $\bx $ and $\by $ is defined to be
$
\langle \bx ,\by\rangle =-x_0y_0+\sum _{i=1}^n x_iy_i.
$
We call $(\R^{n+1} ,\langle ,\rangle )$  the {\it $(n+1)$-dimensional Minkowski space-time} (or briefly, the {\it Lorentz-Minkowski
$(n+1)$-space}). We write $\R^{n+1}_1$ instead of $(\R^{n+1} ,\langle
,\rangle )$. We say that a non-zero vector $\bx \in \R^{n+1}_1$ is
{\it spacelike, lightlike or timelike} if $\langle \bx,\bx \rangle
>0,$ $\langle \bx,\bx \rangle =0$ or  $\langle \bx,\bx \rangle <0$
respectively. The norm of the vector $\bx \in \R^{n+1}_1$ is defined
to be $\|\bx\|=\sqrt{|\langle \bx,\bx\rangle |}.$ We have the canonical
projection $\pi :\R^{n+1}_1\lon \R^n$ defined by $\pi (x_0,x_1,\dots
,x_n)=(x_1,\dots ,x_n).$ Here we identify $\{\bo\}\times \R^n$ with
$\R^n$ and it is considered as the Euclidean $n$-space whose scalar
product is induced from the pseudo scalar product $\langle ,\rangle
.$ For a vector $\bv \in \R^{n+1}_1$ and a real number $c,$ we
define a {\it hyperplane with pseudo normal\/} $\bv$ by
$$
HP(\bv ,c)=\{\bx\in \R^{n+1}_1\ |\ \langle \bx ,\bv \rangle = c\ \}.
$$
We call $HP(\bv ,c)$ a {\it spacelike hyperplane\/}, a {\it timelike hyperplane\/}
or a {\it lightlike hyperplane\/}  if $\bv $ is timelike, spacelike or lightlike respectively.
\par
We now define {\it Hyperbolic $n$-space} by
$$
H^n_+(-1)=\{\bx\in \R^{n+1}_1 | \langle \bx ,\bx\rangle =-1, x_0>0
\}
$$
and {\it de Sitter $n$-space} by
$$
S^n_1=\{\bx\in \R^{n+1}_1 | \langle \bx ,\bx\rangle =1\ \}.
$$
We define
$$
LC(\bm{\lambda})=\{\bx=(x_0,x_1,\dots ,x_n)\in \R^{n+1}_1  \ |\
\langle\bx-\bm{\lambda} ,\bx -\bm{\lambda}\rangle =0\}
$$
and we call it {\it the lightcone} with the vertex $\bm{\lambda}\in \R^{n+1}_1.$
We write $LC^*=LC(\bm{0})\setminus \{\bm{0}\}$, which is called an {\it open lightcone} at the origin.
\par
For any $\bx_1,\bx_2,\dots ,\bx_n \in \R^{n+1}_1,$
we define a vector $\bx_1\wedge\bx_2\wedge\dots \wedge\bx_n$
by
\[
\bx_1\wedge\bx_2\wedge\dots \wedge\bx_n=
\left|
\begin{array}{cccc}
-\be_{0}&\be_{1}&\cdots &\be_{n}\\
x^1_{0}&x^1_{1}&\cdots &x^1_{n}\\
x^2_{0}&x^2_{1}&\cdots &x^2_{n}\\
\vdots &\vdots &\cdots &\vdots \\
x^n_{0}&x^n_{1}&\cdots &x^n_{n}
\end{array}
\right| ,
\]
where $\be_{0},\be_{1},\dots ,\be_{n}$ is the canonical basis of $\R^{n+1}_1$
and $\bx_i=(x_0^i,x_1^i,\dots ,x_n^i).$
We can easily check that
$
\langle \bx,\bx_1\wedge\bx_2\wedge\dots \wedge\bx_n\rangle ={\rm det}(\bx,\bx_1,\dots ,\bx_n),
$
so that
$\bx_1\wedge\bx_2\wedge\dots \wedge\bx_n$ is pseudo orthogonal to  any $\bx _i$ $(i=1,\dots ,n).$
%%%%%%%%%%%%%%%%%%%%%%%%%%%%%%%%%%%%%%%%%%%%%%%%%%%%
\section{World hyper-sheets in the Minkowski space-time}
%%%%%%%%%%%%%%%%%%%%%%%%%%%%%%%%%%%%%%%%%%%%%%%%%%%
\par
We introduce the basic geometrical framework for the 
study of world hyper-sheets in the $(n+1)$-dimensional Minkowski space-time.
Let $\R^{n+1}_1$ be a time-oriented 
space (cf., \cite{Oneil}). We choose $\be _0=(1,0,\dots ,0)$ as the future timelike
vector field. 
In the theory of relativity, we do not have the notion of 
time constant, so that everything that is moving depends on the time.
Therefore, we  consider world sheets.
Although we have the notion of world sheets with general codimension, we stick to the case when
the codimension one, that is world hyper-sheets in the Minkowski space-time.
The world sheet is defined to be a timelike submanifold foliated by
codimension one spacelike submanifolds.
Here, we only consider the local situation, so that we considered a one-parameter family of spacelike submanifolds.
Let $\bX :U\times I\lon \R^{n+1}_1$ be a timelike embedding from an open subset $U\subset \R^{n-1}$ and an open interval $I$.
 We write 
$W=\bX(U\times I)
$ and identify $W$ and $U\times I$ through the embedding $\bX.$
The embedding $\bX$ is said to be {\it timelike} if the tangent space $T_p  W$
of $W$ is a timelike hyperplane at any point $p\in W$.
We write that $\mathcal{S}_t=\bX(U\times\{t\})$ for each $t\in I.$
We have a foliation of $W$ defined by $\mathcal{S}=\{\mathcal{S}_t\}_{t\in I}$.
We say that $W=\bX(U\times I)$ (or, $\bX$ itself) is a {\it world hyper-sheet}
if $W$ is a time-orientable timelike hypersurface and each $\mathcal{S}_t$ is spacelike.
Here, we say that $\mathcal{S}_t$ is {\it spacelike} if the tangent space $T_p\mathcal{S}_t$ 
consists only spacelike vectors (i.e. spacelike subspace) for
any point $p\in \mathcal{S}_t.$
Each $\mathcal{S}_t$ is called a {\it momentary space} of $W$.
For any $p=\bX(\overline{u},t)\in W\subset \R_1^{n+1},$
we have
\[
T_pW=\langle \bX_{u_1}(\ou,t),\dots ,\bX_{u_{n-1}}(\ou,t),\bX_t(\ou,t)\rangle _\R,
\]
where we write $(\overline{u},t)=(u_1,\dots ,u_{n-1},t)\in U\times I$, $\bX_t=\partial \bX/\partial t$
and $\bX_{u_j}=\partial \bX/\partial u_j.$
We also have
\[
T_p\mathcal{S}_t=\langle \bX_{u_1}(\ou,t),\dots ,\bX_{u_{n-1}}(\ou,t)\rangle _\R.
\]
Since $W$ is time-orientable, there exists a timelike vector field $\bv(\ou,t)$ on $W$ \cite[Lemma 32]{Oneil}.
Moreover, we can choose that $\bv$ is {\it future directed} which means that
$\langle \bm{v}(\ou,t),\bm{e}_0\rangle <0.$
Since ${\rm codim}\, W=1,$ we have ${\rm codim}\, \mathcal{S}_t=2.$
Moreover, $\mathcal{S}_t$ is spacelike, so that we can apply the method developed in \cite{IzuSM}.
We consider the unit normal spacelike vector of $W$ defined by
\[
\bm{n}^S(\ou,t)=\frac{\bX_{u_1}(\ou,t)\wedge\cdots \wedge \bX_{u_{n-1}}(\ou,t)\wedge \bX_t(\ou,t)}{\|\bX_{u_1}(\ou,t)\wedge\cdots \wedge \bX_{u_{n-1}}(\ou,t)\wedge \bX_t(\ou,t)\|}.
\]
For any $t\in I,$ Let $N_p(\mathcal{S}_t)$ be the pseudo-normal space of $\mathcal{S}_t$ at $p=\bX(\ou,t)$ in $\R^{n+1}_1.$
Since $\mathcal{S}_t$ is a codimension one in $W,$ $N_p(\mathcal{S}_t)$ is a two dimensional Lorentz space.
There exists a unique timelike unit vector field $\bm{n}^T(\ou,t)\in N_p(\mathcal{S}_t)\cap T_pW$ such that it is future directed (i.e. $\langle \bm{n}^T(\ou,t),\bm{e}_0\rangle <0$).
We now define a map $\mathbb{LG}^\pm (\mathcal{S}_t):\mathcal{S}_t\lon LC^*$ by
$\mathbb{LG}^\pm (\mathcal{S}_t)(p )=\bm{n}^T(\ou,t)\pm \bm{n}^S(\ou,t),$ where $p=\bm{X}(\ou,t).$
We call each one of $\mathbb{LG}^\pm (\mathcal{S}_t)$ a {\it momentary lightcone Gauss map}.
This map leads us to the notion of curvatures (cf. \cite{Izu14}). We have a linear mapping
$d\mathbb{LG}^\pm (\mathcal{S}_t)_p:T_p\mathcal{S}_t\lon T_{\widetilde{p}}LC^*\subset T_{\widetilde{p}}\R^{n+1}_1,$
where $p =\bm{X}(\ou ,t)$ and $\widetilde{p}=\bm{n}^T(\ou,t)\pm \bm{n}^S(\ou,t)$.
With the identification $T_{\widetilde{p}}\R^{n+1}_1\equiv \R^{n+1}_1\equiv T_p\R^{n+1}_1,$
we have the canonical decomposition $T_p\R^{n+1}=T_p\mathcal{S}_t\oplus N_p(\mathcal{S}_t).$
Let $\Pi^t:T_p\R^{n+1}=T_p\mathcal{S}_t\oplus N_p(\mathcal{S}_t)\lon T_p\mathcal{S}_t$
be the canonical projection. Then we have
linear transformations
\[
S^\pm_\ell(\mathcal{S}_t)_p=-\Pi^t\circ d\mathbb{LG}^\pm (\mathcal{S}_t)_p:T_p\mathcal{S}_t\lon T_p\mathcal{S}_t.
\]
Each one of the above mappings is called a {\it momentary lightcone shape operator} of $\mathcal{S}_t$ at $p=\bm{X}(\ou,t).$
Let $\{\kappa ^\pm_i(\mathcal{S}_t)(p)\}_{i=1}^{n-1}$ be the set of eigenvalues of $S^\pm_\ell(\mathcal{S}_t)_p,$ which are called {\it momentary lightcone principal curvatures} of $\mathcal{S}_t$ at $p=\bm{X}(\ou,t).$
Then {\it momentary lightcone Gauss-Kronecker curvatures} of $\mathcal{S}_t$ at $p=\bm{X}(\ou,t)$ are
defined to be 
\[
K_\ell^\pm(\mathcal{S}_t)(p)=\det S^\pm_\ell(\mathcal{S}_t)_p.
\]
\par
We deduce now the lightcone Weingarten formula. Since $\mathcal{S}_t$
is a spacelike submanifold, we have a Riemannian metric
(the {\it first fundamental form \/}) on $\mathcal{S}_t$ defined by
$ds^2 =\sum _{i=1}^{n-1} g_{ij}du_idu_j$,  where $g_{ij}(\ou,t) =\langle
\bX _{u_i}(\ou,t ),\bX _{u_j}(\ou,t)\rangle$ for any $(\ou,t)\in U\times I.$ 
{\it Lightcone second fundamental invariants\/} are defined to be $h[\pm]
_{ij}(\ou,t)=\langle -(\bn^T \pm\bn^S)
_{u_i}(\ou,t),\bX_{u_j}(\ou,t)\rangle$ for any $(\ou,t)\in U\times I.$
The following lightcone
Weingarten formulae are given as special cases of the formulae in \cite{IzuSM}:
\smallskip

{\rm (a)} $(\bn^T \pm \bn^S)_{u_i}=\langle \bn ^S,\bn
^T_{u_i}\rangle(\bn^T\pm \bn^S)-\sum_{j=1}^{n-1} h_i^j[\pm]\bX
_{u_j}$
\smallskip
\par
{\rm (b)} $ \Pi ^t\circ (\bn^T +\bn^S)_{u_i}=-\sum_{j=1}^{n-1}
h_i^j[\pm ]\bX _{u_j}. $
\smallskip
\par\noindent
Here $\displaystyle{\left(h_i^j[\pm]\right)=\left(h_{ik}[\pm]\right)\left(g^{kj}\right)}$ and
$\displaystyle{\left( g^{kj}\right)=\left(g_{kj}\right)^{-1}}.$
\par\noindent
It follows that the momentary lightcone principal curvatures are the eigenvalues of $\displaystyle{\left(h_i^j[\pm]\right)}.$

%%%%%%%%%%%%%%%%%%%%%%%%%%
\section{Light sheets along momentary spaces}
%%%%%%%%%%%%%%%%%%%%%%%%%%%
 We define a hypersurface
$
\mathbb{LH}^\pm_{\mathcal{S}_{t_0}}:U\times\{t_0\}\times \R\lon \R^{n+1}_1$
by
$$
\mathbb{LH}^\pm_{\mathcal{S}_{t_0}}(p,\mu)=\mathbb{LH}^\pm_{\mathcal{S}_{t_0}}(\ou,t_0,\mu)=\bX(\ou,t_0)+\mu\mathbb{LG}^\pm(\mathcal{S}_{t_0})(\ou,t_0),
$$
where $p=\bX (\ou,t_0).$
We call $\mathbb{LH}^\pm_{\mathcal{S}_{t_0}}$ {\it light sheets\/} along $\mathcal{S}_{t_0}.$
In general, a hypersurface $H\subset \R^{n+1}_1$ is called a {\it lightlike hypersurface} if it is tangent to
a lightcone at any point.
The light sheet along $\mathcal{S}_{t_0}$ is a lightlike hypersurface.
We also define $\mathbb{LH}^\pm_W:U\times I\times \R\lon \R^{n+1}_1\times I$ by
\[
\mathbb{LH}^\pm_W(\ou,t,\mu)=(\mathbb{LH}^\pm_{\mathcal{S}_{t}}(\ou, t,\mu),t),
\]
which is called an {\it unfolded light sheets} of $(W,\mathcal{S}).$
 \par
 We introduce the notion of Lorentz distance-squared
functions on a world hyper-sheet, which is useful for the study of
singularities of light sheets.
We define a family of functions $G: W\times\R_1^{n+1}\lon \R$ on $W=\bX (U\times I)$ 
 by
$$
 G(p,\blambda)=G(\ou,t, \blambda )=\langle \bX (\ou,t )-\blambda ,\bX
(\ou,t)-\blambda\rangle ,
$$
where $p=\bX(\ou,t).$
We call
$G$ a {\it Lorentz distance-squared function\/} on the world hyper-sheet
$(W,\mathcal{S}).$
For any fixed $(t_0,\blambda _0)\in I\times \R_1^{n+1},$ we write $g(\ou)=G_{(t_0,\bm{\lambda} _0)}( \ou)=G(\ou, t_0,\blambda _0)$
and have
the following proposition.
\par
\begin{Pro}
Let $\mathcal{S}_{t_0}$ be a momentary space of $(W,\mathcal{S})$ and
$G: W\times\R_1^{n+1}\to\R$
the Lorentz distance-squared function on $(W,\mathcal{S}).$
Suppose that $p_0=\bX(\ou_0,t_0)\not=\blambda _0.$ Then we have the following$:$
\par
{\rm (1)}
$g(\ou_0)=\partial g/\partial u_i(\ou_0)=0$ $(i=1,\dots ,n-1)$
if and only if
$p_0-\blambda _0 =\mu\mathbb{LG}^\pm(\mathcal{S}_{t_0}) (p_0)$ for some $\mu\in
\R\setminus \{0\}.$
\par
{\rm (2)}
$g(\ou_0)=\partial g/\partial u_i(\ou_0)=
{\rm det}\,{\mathcal H}(g)(\ou_0)=0$ $(i=1,\dots ,n-1)$
if and only if
$$
p_0-\blambda _0=\mu\mathbb{LG}^\pm(\mathcal{S}_{t_0}) (p_0)
$$
for $\mu\in \R\setminus \{0\}$ such that
$-1/\mu$  is one of the non-zero momentary lightcone
principal curvatures
$\{\kappa ^\pm_i(\mathcal{S}_t)(p)\}_{i=1}^{n-1}.$
\par\noindent
Here, ${\rm det}\,{\mathcal H}(g)(\ou_0)$
is the determinant of the Hessian matrix of $g$ at $\ou_0.$.
\end{Pro}
\par
\demo
(1) The condition $g(\ou)=\langle \bX(\ou, t_0)-{\blambda _0},\bX(\ou,t_0)-{\blambda
_0}\rangle =0$
means that
$\bX(\ou,t_0)-{\blambda _0}\in LC^*.$
We can observe that $dg(\ou)=\langle d\bX(\ou,t_0), \bX(\ou,t_0)-{\blambda _0}\rangle =0$
if and only if
$\bX (\ou,t_0)-{\blambda _0}\in N_pM.$
Hence
$g(\ou_0)=dg(\ou_0)=0$
if and only if
$p_0-{\blambda _0}\in N_pM\cap LC^*.$
This is equivalent to the condition that
$p_0-{\blambda _0}=\mu\mathbb{LG}^\pm(\mathcal{S}_{t_0})(p_0)$ for some $\mu\in
\R\setminus \{0\}.$
\par
(2) We can calculate that
\[
\frac{\partial g}{\partial u_i}=2\langle \bX_{u_i},\bX -\blambda _0\rangle
\]
and
\[
\frac{\partial ^2 g}{\partial u_i\partial u_j}
=2\left\{\bX _{u_iu_j},\bX-\blambda _0\rangle +\langle \bX_{u_i},\bX _{u_j}\rangle\right\}.
\]
With the condition $p_0-{\blambda _0}=\mu\mathbb{LG}^\pm(\mathcal{S}_{t_0})(p_0)$, we have
\[
\frac{\partial ^2 g}{\partial u_i\partial u_j}
=2\left\{\langle \bX _{u_iu_j},\mu \mathbb{LG}^\pm(\mathcal{S}_{t_0})(p_0)\rangle +g_{ij}(\ou_0,t_0)\right\}.
\]
Therefore, we have
\[
\left(\frac{\partial ^2 g}{\partial u_i\partial u_j}\right)\left(g^{k\ell}\right)
=\left(2\left\{\mu{h}^i_j[\pm]+\delta ^i_j\right\}\right).
\]
It follows that ${\rm det}{\mathcal H}(g)(p_0)=0$ if and only if 
$-1/\mu$ is an eigenvalue of $({h}^i_j[\pm](p_0)).$
\enD

\par
Inspired by the above result, we define 
\[
\mathbb{LF}^\pm_{\mathcal{S}_{t_0}}=\bigcup _{i=1}^{n-1}\left\{\bX(u,t_0)+\frac{1}{\kappa ^\pm_i(\mathcal{S}_t)( p)}\mathbb{LG}^\pm(\mathcal{S}_{t_0})( p)\ |\ 
u\in U, p=\bX(u,t_0)\ \right\},
\]
which are called {\it lightlike focal sets} of $\mathcal{S}_{t_0}.$
Moreover, {\it unfolded lightcone focal sets} of $(W,\mathcal{S})$ are defined to be
\[
\mathbb{LF}^\pm_{(W,\mathcal{S})}=\bigcup _{t\in I} \mathbb{LF}^\pm_{\mathcal{S}_{t}}\times \{t\} \subset \mathbb{R}^{n+1}_1\times I.
\]
Each one of $\mathbb{LF}^\pm_{(W,\mathcal{S})}$ is the critical value set of $\mathbb{LH}^\pm_W$, respectively.
%%%%%%%%%%%%%%%%%%%%%%%%%%%%%%%%%%%
\section{Contact with lightcones}
%%%%%%%%%%%%%%%%%%%%%%%%%%%%%%%%%
In this section we interpret the results of Proposition 4.1 from the view point of the contact with lightcones.
\par
Firstly,
we consider the relationship between the contact of a one parameter family of submanifolds with a submanifold and 
$S.P$-${\mathcal K}$-equivalence among functions (cf., \cite{Izudoc}).  
Let $ U_i\subset \R^r$, ($i=1,2$) be open sets and $g_i:(U_i\times I, (\ou_i,t_i))\lon (\R^n,\bm{y}_i)$  immersion germs. We define $\overline{g}_i:(U_i\times I, (\ou_i,t_i))\lon (\R^n\times I,(\bm{y}_i,t_i))$
by $\overline{g}_i(\ou,t)=(g_i(\ou),t).$
We write that $(\overline{Y}_i,(\bm{y}_i,t_i))=(\overline{g}_i(U_i\times I),(\bm{y}_i,t_i)).$
Let 
$f_i:(\R^n,\bm{y}_i) \lon (\R,0)$ be submersion germs and write that $(V(f_i),\bm{y}_i)=(f_i^{-1}(0),\bm{y}_i).$
We say that {\it the contact of $\overline{Y}_1$ with the trivial family of $V(f_1)$  
at $(\bm{y}_1,t_1)$} is of the {\it same type in the strict sense} as {\it the contact of 
$\overline{Y}_2$ with the trivial family of $V(f_2)$ at $(\bm{y}_2,t_2)$}
if there is a diffeomorphism germ
$\Phi:(\R^n\times I,(\bm{y}_1,t_1)) \lon (\R^n\times I,(\bm{y}_2,t_2))$ of the form $\Phi (\bm{y},t)=(\phi_1(\bm{y},t),t+(t_2-t_1))$
 such that $\Phi(\overline{Y}_1)=\overline{Y}_2$ and
$\Phi(V(f_1)\times I) = V(f_2)\times I$.
In this case we write 
$SK(\overline{Y}_1,V(f_1)\times I;(\bm{y}_1,t_1)) = SK(\overline{Y}_2,V(f_2)\times I;(\bm{y}_2,t_2))$.
We can show one of the parametric versions of Montaldi's theorem of contact between submanifolds as follows: 
\begin{Pro} 
With the same notations as in the above paragraph, 
$SK(\overline{Y}_1,V(f_1)\times I;(\bm{y}_1,t_1)) = SK(\overline{Y}_2,V(f_2)\times I;(\bm{y}_2,t_2))$ 
if and only if 
$f_1 \circ g_1$ and $f_2 \circ g_2$ are $S.P$-${\mathcal K}$-equivalent
 {\rm (}i.e. there exists a diffeomorphism germ $\Psi :(U_1\times I,(\ou_1,t_1))\lon (U_2\times I,(\ou_2,t_2))$ of the form $\Psi (\ou,t)=(\psi_1 (\ou,t),t-(t_2-t_1))$ and a function germ $\lambda :(U_1\times I,(\ou_1,t_1))\lon \R$ with 
 $\lambda (\ou_1,t_1)\not= 0$ such that $(f_2\circ g_2)\circ \Phi (\ou,t) =\lambda (\ou,t)f_1\circ g_1(\ou,t)${\rm ).}
\end{Pro}
Since the proof of Proposition 5.1 is given by the arguments just along the line of the proof of the original theorem in \cite{mont1},
we omit the proof here.
\par
We now consider a function 
$
{\mathfrak g}_{\bm{\lambda}}:\R^{n+1}_1\lon {\mathbb R}
$
defined by ${\mathfrak g}_{\bm{\lambda}}(\bx )=\langle \bx-\bm{\lambda} ,\bx -\bm{\lambda}\rangle,$
where $\bm{\lambda}\in \R^{n+1}_1\setminus W.$
For any $\bm{\lambda}_0\in \R^{n+1}_1$, we have a lightcone $\mathfrak{g}_{\bm{\lambda}_0}^{-1}(0)=LC(\bm{\lambda}_0).$
Moreover, we consider the lightlike vector $\bm{\lambda}_0=\mathbb{LH}^\pm_{\mathcal{S}_{t_0}} (p_0,\mu_0),$
where $p_0=\bX(\ou_0,t_0).$
Then we have 
$$
{\mathfrak{g}}_{\bm{\lambda}_0}\circ\bX (\ou_0,t_0)=G((u_0,t_0), \mathbb{LH}^\pm_{\mathcal{S}_{t_0}} (p_0,\mu_0))=0.
$$
By Proposition 4.1, we also have relations that
$$
\frac{\partial {\mathfrak{g}}_{\bm{\lambda}_0}\circ\bX }{\partial u_i}(\ou_0,t_0)=\frac{\partial G}{\partial u_i}((u_0,t_0), \mathbb{LH}^\pm_{\mathcal{S}_{t_0}} (p_0,\mu_0))=0.
$$
for $i=1,\dots ,n-1.$
This means that the lightcone ${\mathfrak{g}}_{\bm{\lambda}_0}^{-1}(0)=LC(\bm{\lambda}_0)$ is tangent to 
$\mathcal{S}_{t_0}=\bX (U\times \{t_0\})$ at $p_0=\bX (\ou_0,t_0).$
The lightcone $LC(\bm{\lambda}_0)$ is said to be a {\it tangent lightcone} of $\mathcal{S}_{t_0}=\bX (U\times \{t_0\})$ at 
$p_0=\bX (\ou_0,t_0)$, which we write 
$TLC(\mathcal{S}_{t_0},\bm{\lambda}_0),$ where $\bm{\lambda}_0=\mathbb{LH}^\pm_{\mathcal{S}_{t_0}} (p_0,\mu_0).$
Then we have the following simple lemma.
\begin{Lem}
Let $\bX :U\times I\lon \R^{n+1}_1$ be a world hyper-sheet. Consider two points $p_i=\bX(\ou_i,t_0)$, $(i=1,2).$
Then \[
\mathbb{LH}^\pm_{\mathcal{S}_{t_0}} (p_1,\mu_1)=\mathbb{LH}^\pm_{\mathcal{S}_{t_0}}(p_2,\mu_2)
\] if and only if 
$$TLC(\mathcal{S}_{t_0},\mathbb{LH}^\pm_{\mathcal{S}_{t_0}} (p_1,\mu_1))=TLP(\mathcal{S}_{t_0},\mathbb{LH}^\pm_{\mathcal{S}_{t_0}}(p_2,\mu_2)).$$
\end{Lem}
\par
Eventually, we have tools for the study of the contact between momentary spaces and lightcones.
Since we have $g_{\bm{\lambda}}(\ou,t)=\mathfrak{g}_{\bm{\lambda}}\circ\bX(\ou,t),$ we have the following proposition as a
corollary of Proposition 5.1.
\begin{Pro}
Let $\bX_i : (U\times I,(\ou_i,t_0)) \lon (\R^{n+1}_1,p_i)$, $(i=1,2)$, be  world hypersheet germs and $\bm{\lambda}_i=\mathbb{LH}^\pm_{\mathcal{S}_{t_0}} (p_i,\mu_i))$
and $W_i=\bX_i(U\times I).$
Then the following conditions are equivalent:
\par\noindent
{\rm (1)} $SK(\overline{W}_1, TLC(\mathcal{S}_{t_0},\bm{\lambda}_1)\times I;(p_1,t_0))=
SK(\overline{W}_2, TLC(\mathcal{S}_{t_0},\bm{\lambda}_2)\times I;(p_2,t_0)),$ 
\par\noindent
{\rm (2)} $g_{1,\bm{\lambda}_1}$ and $g_{2,\bm{\lambda}_2}$ are $S.P$-$\mathcal{K}$-equivalent.\\
Here, $g_{i,\bm{\lambda}_i}(\ou,t)=\langle \bX_i(\ou,t)-\bm{\lambda}_i,\bX_i(\ou,t)-\bm{\lambda}_i\rangle ,$
$(i=1,2)$.
\end{Pro}

%%%%%%%%%%%%%%%%%%%%%%%%%%%%%%%%%
\section{Graph-like wave fronts}
%%%%%%%%%%%%%%%%%%%%%%%%%%%%%%%%%%
In this section we briefly review the theory of
graph-like Legendrian unfoldings.
Graph-like Legendrian unfoldings belong to a special class of big Legendrian 
submanifolds (for detail, see \cite{Izumiya93,Izumiya-Takahashi,Izumiya-Takahashi2,Izumiya-Takahashi3,Zakalyukin95}).
Recently there appeared a survey article \cite{Graph-like} on the theory of graph-like Legendrian unfoldings.
Let ${\mathcal F} :(\R^k\times (\R^m\times\R),0)\to (\R,0)$ be a function germ.
We say that ${\mathcal F}$ is a {\it graph-like Morse family of hypersurfaces}
if $
(\mathcal{F}, d_q\mathcal{F}):(\R^k\times(\R^m\times \R),0)\to (\R\times \R^k,0)$
is a non-singular and 
$(\partial {\mathcal F}/\partial t)(0)\not= 0,$
where
$$
d_q\mathcal{F}(q,x,t)=\left(\frac{\partial \mathcal{F}}{\partial q_1}(q,x,t), \dots ,
\frac{\partial \mathcal{F}}{\partial q_k}(q,x,t)\right).
$$
Moreover, we say that ${\mathcal F}$ is {\it non-degenerate} if
$(\mathcal{F}, d_q\mathcal{F})|_{\R^k\times(\R^m\times \{0\})}$ is non-singular.
For a graph-like Morse family of hypersurfaces $\mathcal{F},$ $\Sigma _*(\mathcal{F})=(\mathcal{F}, d_q\mathcal{F})^{-1}(0)$ is a
smooth $m$-dimensional submanifold germ of $(\R^k\times(\R^m\times \R),0).$
We now consider the space of $1$-jets $J^1(\R^m,\R)$ with the canonical coordinates $(x_1,\dots ,x_m,t,p_1,\dots ,p_m)$ such
that the canonical contact form is $\theta =dt-\sum_{i=1}^m p_idx_i.$
We define a mapping $\Pi:J^1(\R^m,\R)\lon T^*\R^m$ by
$\Pi(x,t,p)=(x,p),$ where $(x,t,p)=(x_1,\dots, x_m,t,p_1,\dots ,p_m).$ 
Here, $T^*\R^m$ is a symplectic manifold with the
canonical symplectic structure $\omega=\sum_{i=1}^m dp_i\wedge dx_i$ (cf. \cite{Arnold1}).
We define a mapping $\mathscr{L}_{\mathcal{F}}:(\Sigma_{*}(\mathcal{F}),0) \to J^1(\R^m,\R)$ by 
\[
\mathscr{L}_{\mathcal F}(q,x,t)=\left(x,t,-\frac{\displaystyle \frac{\partial \mathcal{F}}{\displaystyle\partial x_1}(q,x,t)}{\frac{\displaystyle \partial \mathcal{F}}{\displaystyle\partial t}(q,x,t)},\dots , -\frac{\displaystyle\frac{\partial \mathcal{F}}{\partial x_m}(q,x,t)}{\frac{\displaystyle\partial \mathcal{F}}{\displaystyle\partial t}(q,x,t)},\right).
\]
It is easy to show that $\mathscr{L}_{\mathcal F}(\Sigma _*({\mathcal F}))$ is a Legendrian submanifold germ (cf., \cite{Arnold1}), which is called a {\it graph-like Legendrian unfolding germ.} 
We call $\overline{\pi}|_{\mathscr{L}_{\mathcal F}(\Sigma _*({\mathcal F}))}:\mathscr{L}_{\mathcal F}(\Sigma _*({\mathcal F}))\lon \R^m\times \R$  a {\it graph-like Legendrian map germ}, where
$\overline{\pi}:J^1(\R^m,\R)\lon \R^m\times\R$ is the canonical projection.
We also call $W(\mathscr{L}_{\mathcal F}(\Sigma _*({\mathcal F})))=\overline{\pi}(\mathscr{L}_{\mathcal F}(\Sigma _*({\mathcal F})))$ a {\it graph-like wave front}
of $\mathscr{L}_{\mathcal F}(\Sigma _*({\mathcal F})).$
We say that ${\mathcal F}$ is a {\it graph-like generating family} of
$\mathscr{L}_{\mathcal F}(\Sigma _*({\mathcal F})).$
Moreover, we call $W_t(\mathscr{L}_{\mathcal F}(\Sigma _*({\mathcal F})))=\pi _1(\pi _2^{-1}(t)\cap W(\mathscr{L}_{\mathcal F}(\Sigma _*({\mathcal F})))$
a {\it momentary front} for each $t\in (\R,0),$ where $\pi _1:\R^m\times \R\lon \R^m$ and $\pi _2:\R^m\times \R\lon \R$ are the canonical projections.
The {\it discriminant set of the family} $\{W_t(\mathscr{L}_{\mathcal F}(\Sigma _*({\mathcal F})))\}_{t\in (\R,0)}$ is defined by the union of the caustic $C_{\mathscr{L}_{\mathcal F}(\Sigma _*({\mathcal F}))}=\pi _1( \Sigma (W (\mathscr{L}_{\mathcal F}(\Sigma _*({\mathcal F}))))$  and the Maxwell stratified set $M_{\mathscr{L}_{\mathcal F}(\Sigma _*({\mathcal F}))}=\pi _1(SI_{W(\mathscr{L}_{\mathcal F}(\Sigma _*({\mathcal F})))}),$ where $\Sigma (W (\mathscr{L}_{\mathcal F}(\Sigma _*({\mathcal F})))$ is the critical value set of $\overline{\pi}|_{\mathscr{L}_{\mathcal F}(\Sigma _*({\mathcal F}))}$
and $SI_{W(\mathscr{L}_{\mathcal F}(\Sigma _*({\mathcal F})))}$ is the closure of the self intersection set of $W(\mathscr{L}_{\mathcal F}(\Sigma _*({\mathcal F}))).$
\par
 We now define equivalence relations among graph-like Legendrian unfoldings.
Let ${\mathcal F} :(\R^k\times (\R^m\times\R),0)\to (\R,0)$ and ${\mathcal G} :(\R^{k}\times (\R^m\times\R),0)\to (\R,0)$ be graph-like Morse families of hypersurfaces.
We say that $\mathscr{L}_{\mathcal F}(\Sigma _*({\mathcal F}))$ and $\mathscr{L}_{\mathcal G}(\Sigma _*({\mathcal G}))$
are {\it Legendrian equivalent} if there exist a diffeomorphism germ $\Phi:(\R^m\times \R,\overline{\pi}(p ))\lon (\R^m\times \R,\overline{\pi}(p' ))$ and
a contact diffeomorphism germ $\widehat{\Phi}:(J^1(\R^m,\R),p)\lon (J^1(\R^m,\R),p')$ such that $\overline{\pi}\circ\widehat{\Phi}=\Phi\circ \overline{\pi}$
and $\widehat{\Phi}(\mathscr{L}_{\mathcal F}(\Sigma _*({\mathcal F})))=(\mathscr{L}_{\mathcal G}(\Sigma _*({\mathcal G}))),$
where $p=\mathscr{L}_{\mathcal F}(0)$ and $p'=\mathscr{L}_{\mathcal G}(0).$
We also say that $\mathscr{L}_{\mathcal F}(\Sigma _*({\mathcal F}))$ and $\mathscr{L}_{\mathcal G}(\Sigma _*({\mathcal G}))$
are {\it $S.P^+$-Legendrian equivalent} if these are Legendrian equivalent by a diffeomorphism germ $\Phi:(\R^m\times \R,\overline{\pi}(p ))\lon (\R^m\times \R,\overline{\pi}(p' ))$ 
of the form $\Phi (x,t)=(\phi _1(x),t+\alpha (x))$ and 
a contact diffeomorphism germ $\widehat{\Phi}:(J^1(\R^m,\R),p)\lon (J^1(\R^m,\R),p')$ with $\overline{\pi}\circ\widehat{\Phi}=\Phi\circ \overline{\pi}.$
Moreover, graph-like wave fronts $W(\mathscr{L}_{\mathcal F}(\Sigma _*({\mathcal F})))$ and $W(\mathscr{L}_{\mathcal G}(\Sigma _*({\mathcal G})))$
are {\it $S.P^+$-diffeomorphic} if there exists a diffeomorphism germ $\Phi:(\R^m\times \R,\overline{\pi}(p ))\lon (\R^m\times \R,\overline{\pi}(p' ))$ 
of the form $\Phi (x,t)=(\phi _1(x),t+\alpha (x))$ such that $\Phi(W(\mathscr{L}_{\mathcal F}(\Sigma _*({\mathcal F}))))=W(\mathscr{L}_{\mathcal G}(\Sigma _*({\mathcal G})))$
as set germs.
By definition, if $\mathscr{L}_{\mathcal F}(\Sigma _*({\mathcal F}))$ and $\mathscr{L}_{\mathcal G}(\Sigma _*({\mathcal G}))$
are $S.P^+$-Legendrian equivalent, then $W(\mathscr{L}_{\mathcal F}(\Sigma _*({\mathcal F})))$ and $W(\mathscr{L}_{\mathcal G}(\Sigma _*({\mathcal G})))$
are $S.P^+$-diffeomorphic. The converse assertion holds generically \cite{Graph-like,GeomLag14}.
\begin{Pro}[\cite{GeomLag14}] With the same notations as those of the above, we suppose that the sets of critical points of $\overline{\pi}|_{\mathscr{L}_{\mathcal F}(\Sigma _*({\mathcal F}))}, \overline{\pi}|_{\mathscr{L}_{\mathcal G}(\Sigma _*({\mathcal G}))}$
 are nowhere dense respectively.
 Then $\mathscr{L}_{\mathcal F}(\Sigma _*({\mathcal F}))$ and $\mathscr{L}_{\mathcal G}(\Sigma _*({\mathcal G}))$are $S.P^+$-Legendrian equivalent  if and only if $W(\mathscr{L}_{\mathcal F}(\Sigma _*({\mathcal F})))$ and $W(\mathscr{L}_{\mathcal G}(\Sigma _*({\mathcal G})))$
 are $S.P^+$-diffeomorphic.
\end{Pro}
\par
We remark that if $W(\mathscr{L}_{\mathcal F}(\Sigma _*({\mathcal F})))$ and $W(\mathscr{L}_{\mathcal G}(\Sigma _*({\mathcal G})))$
 are $S.P^+$-diffeomorphic by  a diffeomorphism germ $\Phi:(\R^m\times \R,\overline{\pi}(p ))\lon (\R^m\times \R,\overline{\pi}(p' ))$,
 then $$\Phi ( C_{\mathscr{L}_{\mathcal F}(\Sigma _*({\mathcal F}))}\cup M_{\mathscr{L}_{\mathcal F}(\Sigma _*({\mathcal F}))})=
 C_{\mathscr{L}_{\mathcal G}(\Sigma _*({\mathcal G}))}\cup M_{\mathscr{L}_{\mathcal G}(\Sigma _*({\mathcal G}))}.
 $$

\par
For a graph-like Morse family hypersurfaces $\mathcal{F}:(\R^k\times (\R^m\times\R),0)\to (\R,0),$ by the implicit function theorem, there exist function germs
$F:(\R^k\times \R^m,0)\to (\R,0)$ and $\lambda :(\R^k\times (\R^m\times \R),0)\lon \R$ with $\lambda (0)\not= 0$
such that ${\mathcal F}(q,x,t)
=\lambda(q,x,t)( F(q,x)-t).$
We have shown in \cite{Graph-like} that
$\mathcal{F}$ is a graph-like Morse family of hypersurfaces if and only if
$F$ is a Morse family of functions.
Here we say that $F:(\R^k\times\R^m,0)\lon (\R,0)$ is a {\it Morse family of functions} if
\[
dF_q=\left(\frac{\partial F}{\partial q_1},\dots, \frac{\partial F}{\partial q_k}\right):(\R^k\times\R^m,0)\lon \R^k
\]
is non-singular.
We consider a graph-like Morse family of hypersurfaces $\mathcal{F}(q,x,t)=\lambda (q,x,t)(F(q,x)-t).$ 
In this case, 
\[
\Sigma _*(\mathcal{F})=\{(q,x,F(q,x))\in (\R^k\times (\R^m\times\R),0)\ |\ (q,x)\in C(F)\},
\]
where
\[
C(F)=\left\{ (q,x)\in (\R^k\times\R^m,0)\ \Bigm|\ \frac{\partial F}{\partial q_1}(q,x)=\cdots =\frac{\partial F}{\partial q_k}(q,x)=0\ \right\}.
\]
Moreover, we define a map germ
$
L(F):(C(F),0)\lon T^*\R^m
$
by 
\[
L(F)(q,x)=\left(x,\frac{\partial F}{\partial x_1}(q,x),\dots ,\frac{\partial F}{\partial x_m}(q,x)\right)
\]
It is known that $L(F)(C(F))$ is a Lagrangian submanifold germ (cf., \cite{Arnold1})
for the canonical symplectic structure.
In this case $F$ is said to be a {\it generating family} of the Lagrangian submanifold germ $L(F)(C(F)).$
We remark that $\Pi (\mathscr{L}_{\mathcal F}(\Sigma _*({\mathcal F})))=L(F)(C(F))$ and the graph-like wave front
$W(\mathscr{L}_{\mathcal F}(\Sigma _*({\mathcal F})))$ is the graph of $F|C(F).$
Here we call $\pi|_{L(F)(C(F))}:L(F)(C(F))\lon \R^m$ a {\it Lagrangian map germ}, where
$\pi :T^*\R^m\lon \R^m$ is the canonical projection.
Then the set of critical values of $\pi|_{L(F)(C(F))}$ is called a {\it caustic} of $L(F)(C(F))=\Pi(\mathscr{L}_{\mathcal{F}}(\Sigma _*(\mathcal{F})))$ in the theory of Lagrangian singularities,
which is denoted by $C_{L(F)(C(F))}.$
By definition, we have  $C_{L(F)(C(F))}=C_{\mathscr{L}_{\mathcal F}(\Sigma _*({\mathcal F}))}.$
\par
Let ${\mathcal F}, \mathcal{G} :(\R^k\times (\R^m\times\R),0)\to (\R,0)$ be graph-like Morse families of hypersurfaces.
We say that $\Pi(\mathscr{L}_{\mathcal{F}}(\Sigma _*(\mathcal{F})))$ and $\Pi(\mathscr{L}_{\mathcal{G}}(\Sigma _*(\mathcal{G})))$
are {\it Lagrangian equivalent} if there exist a diffeomorphism germ $\Psi:(\R^m, \pi\circ\Pi(p ))\lon (\R^m,\pi\circ\Pi(p' ))$ and
a symplectic diffeomorphism germ $\widehat{\Psi}:(T^*\R^m,\Pi( p))\lon (T^*\R^m,\Pi( p'))$ such that $\pi\circ\widehat{\Psi}=\Psi\circ \pi$
and $\widehat{\Psi}(\Pi(\mathscr{L}_{\mathcal{F}}(\Sigma _*(\mathcal{F}))))=\Pi(\mathscr{L}_{\mathcal{G}}(\Sigma _*(\mathcal{G}))),$
where $p=\mathscr{L}_{\mathcal F}(0)$ and $p'=\mathscr{L}_{\mathcal G}(0).$
By definition, if $\Pi(\mathscr{L}_{\mathcal{F}}(\Sigma _*(\mathcal{F})))$ and $\Pi(\mathscr{L}_{\mathcal{G}}(\Sigma _*(\mathcal{G})))$ are  Lagrangian equivalent, then the caustics $C_{\mathscr{L}_{\mathcal F}(\Sigma _*({\mathcal F}))}$ and $C_{\mathscr{L}_{\mathcal G}(\Sigma _*({\mathcal G}))}$ are diffeomorphic as set germs.
The converse assertion, however, does not hold (cf. \cite{GeomLag14}).
Recently, we have shown the following theorem (cf. \cite{Izumiya-Takahashi2,Graph-like, GeomLag14})
\begin{Th} With the same notations as the above, $\Pi(\mathscr{L}_{\mathcal{F}}(\Sigma _*(\mathcal{F})))$ and $\Pi(\mathscr{L}_{\mathcal{G}}(\Sigma _*(\mathcal{G})))$
are  Lagrangian equivalent if and only if 
$\mathscr{L}_{\mathcal F}(\Sigma _*({\mathcal F}))$ and $\mathscr{L}_{\mathcal G}(\Sigma _*({\mathcal G}))$
are $S.P^+$-Legendrian equivalent.
\end{Th}
We have the following corollary of Proposition 6.1 and Theorem 6.2.
\begin{Co} Suppose that the sets of critical points of $\overline{\pi}|_{\mathscr{L}_{\mathcal F}(\Sigma _*({\mathcal F}))}, \overline{\pi}|_{\mathscr{L}_{\mathcal G}(\Sigma _*({\mathcal G}))}$
 are nowhere dense, respectively.
 Then $\Pi(\mathscr{L}_{\mathcal{F}}(\Sigma _*(\mathcal{F})))$ and $\Pi(\mathscr{L}_{\mathcal{G}}(\Sigma _*(\mathcal{G})))$
are  Lagrangian equivalent  if and only if $W(\mathscr{L}_{\mathcal F}(\Sigma _*({\mathcal F})))$ and $W(\mathscr{L}_{\mathcal G}(\Sigma _*({\mathcal G})))$
 are $S.P^+$-diffeomorphic.
\end{Co}
There are the notions of Lagrangian stability of Lagrangian submanifold germs and $S.P^+$-Legendrian stability of graph-like Legendrian unfolding germs, respectively.
Here we do not use the exact definitions of those notions of stability, so that we omit to give the definitions.
For detailed properties of such stabilities, see \cite{Arnold1, Graph-like}.
We have the following corollary of Theorem 6.2.
\begin{Co} The graph-like Legendrian unfolding $\mathscr{L}_{\mathcal F}(\Sigma _*({\mathcal F}))$ is $S.P^+$-Legendrian stable if and only if 
the corresponding Lagrangian submanifold $\Pi(\mathscr{L}_{\mathcal{F}}(\Sigma _*(\mathcal{F})))$ is Lagrangian stable.
\end{Co}
\par
 Let ${\mathcal F}:(\R^k\times (\R^m\times\R),0)\to (\R,0)$ be a graph-like Morse family of hypersurfaces. We define $\overline{f}:(\R^k\times \R,0)\lon (\R,0)$ by
$\overline{f}(q,t)={\mathcal F}(q,0,t).$
For graph-like Morse families of hypersurfaces ${\mathcal F} :(\R^k\times (\R^m\times\R),0)\to (\R,0)$ and ${\mathcal G} :(\R^{k}\times (\R^m\times\R),0)\to (\R,0)$,
we say that $\overline{f}$ and $\overline{g}$ are {\it $S.P$-$\mathcal{K}$-equivalent} if there exist a function germ $\nu :(\R^k\times\R,0)\lon \R$ with $\nu(0)\not= 0$ and
a diffeomorphism germ $\phi :(\R^k\times \R,0)\lon (\R^k\times \R,0)$ of the form $\phi (q,t)=(\phi_1(q,t),t)$ such that $\overline{f}(q,t)=\nu (q,t)\overline{g}(\phi (q,t)).$
Although we do not give the definition of $S.P^+$-Legendrian stability, we give a corresponding notion for graph-like Morse family of hypersurfaces.
We say that ${\mathcal F}$ is an {\it infinitesimally $S.P^+$-$\mathcal{K}$-versal unfolding} of $\overline{f}$ if 
$$
{\cal E}_{k+1}=\left\langle \frac{\partial \overline{f}}{\partial q_1},\dots, \frac{\partial \overline{f}}{\partial q_k},\overline{f} \right\rangle _{{\cal E}_{k+1}}+
\left\langle \frac{\partial \overline{f}}{\partial t} \right\rangle _{\R}+
\left\langle \frac{\partial \mathcal{F}}{\partial x_1}|_{\R^k\times\{0\}\times \R} ,\dots
,\frac{\partial \mathcal{F}}{\partial x_m}|_{\R^k\times\{0\}\times\R} \right\rangle _{\R},
$$
where ${\cal E}_{k+1}$ is the local $\R$-algebra of $C^\infty$-function germs $(\R^k\times\R,0)\lon \R.$
It is known the following theorem in \cite{Izu95, Zakalyukin95}.
\begin{Th}
The graph-like Legendrian unfolding
$\mathscr{L} _\mathcal{F}(\Sigma _*(\mathcal{F}))$ is $S.P^+$-Legendre stable if and only if $\mathcal{F}$ is an infinitesimally $S.P^+$-${\cal K}$-versal unfolding of $\overline{f}.$
\end{Th}

In \cite{Graph-like} we have shown the following theorem.
\begin{Th} Let ${\mathcal F},\mathcal{G} :(\R^k\times (\R^m\times\R),0)\to (\R,0)$  be graph-like Morse families of hypersurfaces such that ${\mathscr{L}_{\mathcal F}(\Sigma _*({\mathcal F}))}, {\mathscr{L}_{\mathcal G}(\Sigma _*({\mathcal G}))}$ are $S.P^+$-Legendrian stable. Then the following conditions are equivalent{\rm :}
\par\noindent
{\rm (1)} $\mathscr{L}_{\mathcal{F}}(\Sigma _*(\mathcal{F}))$
and $\mathscr{L}_{\mathcal{G}}(\Sigma _*(\mathcal{G}))$ are $S.P^+$-Legendrian equivalent,
\par\noindent
{\rm (2)} $\overline{f}$ and $\overline{g}$ are $S.P$-$\mathcal{K}$-equivalent,
\par\noindent
{\rm (3)} $\Pi(\mathscr{L}_{\mathcal{F}}(\Sigma _*(\mathcal{F})))$ and $\Pi(\mathscr{L}_{\mathcal{G}}(\Sigma _*(\mathcal{G})))$ are Lagrangian equivalent,
\par\noindent
{\rm (4)} $W(\mathscr{L}_{\mathcal{F}}(\Sigma _*(\mathcal{F})))$ and $W(\mathscr{L}_{\mathcal{G}}(\Sigma _*(\mathcal{G})))$
are $S.P^+$-diffeomorphic.
\end{Th}
%%%%%%%%%%%%%%%%%%%%%%%%%%%%%%%%%%%%%%%%%%%%%%%%%%%%%%%%%%%%%%%%
\section{Unfolded lightcone focal sets of world hyper-sheets}
%%%%%%%%%%%%%%%%%%%%%%%%%%%%%%%%%%%%%%%%%%%%%%%%%%%%%%%%%%%%%%%%
In this section we investigate unfolded lightcone focal sets of world hyper-sheets as an application of the theory of graph-like Legendrian unfoldings.
Firstly, we show the following proposition.
\begin{Pro}
Let $G:U\times I\times (\R^{n+1}_1\setminus W)\to \R$ be a Lorentz distance-squared function on a world hyper-sheet $(W,\mathcal{S}).$
For any point $(\ou_0,t_0, \blambda_0 )\in \Sigma _*(G),$ $G$ is a non-degenerate graph-like Morse family of hypersurfaces around
$(\ou_0,t_0,\blambda_0 ).$
\end{Pro}
\demo
We write that
$$
\bX(\ou,t)=(X_0(\ou,t),X_1(\ou,t),\dots ,X_n(\ou,t))\ {\rm and}\
\blambda =(\lambda _0,\lambda _1,\dots ,\lambda _n).
$$
By definition, we have
$$
G(\ou,t,\blambda )=-(X_0(\ou,t)-\lambda _0)^2+(X_1(\ou,t)-\lambda
_1)^2+\cdots 
+(X_n(\ou,t)-\lambda _n)^2.
$$
We now prove that the mapping $$
\Delta^*G(\ou,t_0,\blambda)=\left(G(\ou,t_0,\blambda), \frac{\partial G}{\partial u_1}(\ou,t_0,\blambda),\dots ,\frac{\partial G}{\partial
u_{n-1}}(\ou,t_0,\blambda)\right)
$$
is non-singular at $(\ou_0,t_0,\blambda_0 )\in \Sigma _*(G).$
Indeed, the Jacobian matrix of $\Delta ^*G|_{U\times\{t_0\}\times \R^{n+1}_1}$ is given by
\newfont{\bg}{cmr10 scaled\magstep5}
\newcommand{\bigA}{\smash{\lower1.0ex\hbox{\bg A}}}
\[
\left(
\begin{array}{ccccc}
 &
2(X_0-\lambda _0) & -2(X_1-\lambda _1) &  \cdots  &
-2(X_n-\lambda _n) \\
\bigA & 2X_{0u_1} &
-2X_{1u_1} & \cdots  &-2X_{nu_1}\\
&\vdots & \vdots & \ddots & \vdots \\
& 2X_{0u_{n-1}} &
-2X_{1u_{n-1}} &\cdots  &-2X_{nu_{n-1}}
\end{array}
\right) ,
\]
where $\bigA$ is the following matrix:
\begin{eqnarray*}
\left(\!\!
\begin{array}{ccc}
2\langle \bX -\blambda ,\bX_{u_1}\rangle & \!\! \cdots\!\! & 2\langle \bX -\blambda ,\bX _{u_{n-1}}\rangle \\
2(\langle \bX_{u_1},\bX_{u_1}\rangle +\langle \bX-\blambda ,\bX_{u_1u_1}\rangle) & \!\! \cdots\!\!  &
2(\langle \bX_{u_1},\bX _{u_{n-1}}\rangle +\langle \bX-\blambda ,\bX_{u_1u_{n-1}}\rangle )\\
\vdots & \!\!\ddots\!\! & \vdots \\
2(\langle \bX_{u_{n-1}},\bX_{u_1}\rangle +\langle \bX-\blambda ,\bX_{u_{n-1}u_1}\rangle) &\!\! \cdots\!\! &
2(\langle \bX_{u_{n-1}},\bX _{u_{n-1}}\rangle +\langle \bX-\blambda ,\bX_{u_{n-1}u_{n-1}}\rangle )
\end{array}
\!\!
\right) .
\end{eqnarray*}

Since $\bX $ is an immersion, the rank of the matrix
\begin{eqnarray*}
\left(
\begin{array}{cccc}
2X_{0u_1} &
-2X_{1u_1} & \cdots  &-2X_{nu_1}\\
\vdots & \vdots & \ddots & \vdots \\
 2X_{0u_{n-1}} &
-2X_{1u_{n-1}} &\cdots  &-2X_{nu_{n-1}}
 \end{array}
\right) .
\end{eqnarray*}
is equal to $n-1.$
Since $\bX -\blambda $ is lightlike and $T_p\mathcal{S}_{t_0}$ is spacelike, then  $\{\bX -\blambda ,\bX_{u_1},\dots , \bX_{u_{n-1}}\}$
is linearly independent at $(u_0,t_0,\blambda_0 )\in \Sigma _*(G).$
This means that the rank of
the matrix
\begin{eqnarray*}
\left(
\begin{array}{cccc}
2(X_0-\lambda _0) & -2(X_1-\lambda _1) &  \cdots  &
-2(X_n-\lambda _n) \\
 2X_{0u_1} &
-2X_{1u_1} & \cdots  &-2X_{nu_1}\\
\vdots & \vdots & \ddots & \vdots \\
2X_{0u_{n-1}} &
-2X_{1u_{n-1}} &\cdots  &-2X_{nu_{n-1}}
 \end{array}
\right) 
\end{eqnarray*}
is equal to $n.$
Therefore the Jacobi matrix of $\Delta^*G|_{U\times\{t_0\}\times \R^{n+1}_1}$
is non-singular at $(u_0,t_0,\blambda_0 )\in \Sigma _*(G).$
\par
On the other hand, we have
\[
\frac{\partial G}{\partial t}(\ou,t,\blambda)=2\langle \bX_t(\ou,t),\bX(\ou,t)-\blambda\rangle.
\]
For any $(\ou_0,t_0,\blambda _0)\in \Sigma _*(G),$
there exists $\mu\not= 0$ such that $\blambda _0=\bX(\ou_0,t_0)+\mu\mathbb{LG}^\pm(\mathcal{S}_{t_0}) (\ou_0,t_0).$
Since $\bn^S(\ou_0,t_0)$ is the unit normal vector of $W,$ we have
\[
\langle \bX_t(\ou_0,t_0),\mathbb{LG}^\pm(\mathcal{S}_{t_0}) (\ou_0,t_0)\rangle =\langle \bX_t(\ou_0,t_0),\bn^T(\ou_0,t_0)\rangle.
\]
Moreover, $\{\bX_t(\ou_0,t_0), \bX_{u_1}(\ou_0,t_0), \dots ,\bX_{u_{n-1}}(\ou_0,t_0)\}$ is a basis of $T_pW$
and $\bn^T(\ou_0,t_0)\in N_p(\mathcal{S}_{t_0})\cap T_pW,$
where $p=\bX(\ou_0,t_0).$ It follows that $\langle \bX_t(\ou_0,t_0),\bn^T(\ou_0,t_0)\rangle\not= 0.$
Then we have
\[
\frac{\partial G}{\partial t}(\ou_0,t_0,\blambda_0)=2\langle \bX_t(\ou_0,t_0),-\mu\mathbb{LG}^\pm(\mathcal{S}_{t_0}) (\ou_0,t_0) \rangle
=-\mu\langle\bX_t(\ou_0,t_0),\bn^T(\ou_0,t_0)\rangle\not= 0.
\]
This completes the proof.
\enD
\par
By Proposition 4.1, we have
\[
\Sigma _*(G)=\{(\ou,t,\mathbb{LH}^\pm_{\mathcal{S}_t}(p,\mu))\in U\times I\times \R^{n+1}_1\ |\ p=\bm{X}(\ou, t), \mu \in \R\setminus \{0\}\}. 
\]
We define a map $\mathscr{L}_{G}:\Sigma _*(G)\lon J^1(\R^{n+1}_1,I)$ by
\[
\mathscr{L}_{G}(\ou,t,\mathbb{LH}^\pm_{\mathcal{S}_t}(p,\mu))=\left(\mathbb{LH}^\pm_{\mathcal{S}_t}(p,\mu), t, \frac{2}{\langle \bX_t(\ou,t),\bn^T(\ou,t)\rangle}\overline{\mathbb{LG}^\pm(\mathcal{S}_t)(\ou,t)}\right),
\]
where we define $\overline{\bm{x}}=(-x_0,x_1,\dots 
x_n)$ for $\bm{x}=(x_0,x_1,\dots ,x_n)\in \R^{n+1}_1.$
By the construction of the graph-like Legendrian unfolding from a graph-like Morse family of hypersurfaces, $\mathscr{L}_{G}(\Sigma _*(G))$ is a graph-like Legendrian unfolding
in $J^1(\R^{n+1}_1,I).$
Therefore, the graph-like wave front is 
\[
W(\mathscr{L}_{G}(\Sigma _*(G)))=\{(\mathbb{LH}^\pm_{\mathcal{S}_t}(p,\mu), t)\in \R^{n+1}_1\times I\ |\  p=\bX(\ou,t), (\ou,t)\in U\times I, \mu\in \R\setminus \{0\}\}.
\]
This means that $W(\mathscr{L}_{G}(\Sigma _*(G)))=\mathbb{LH}^+_W(U\times I\times (\R\setminus \{0\}))\cup \mathbb{LH}^-_W(U\times I\times (\R\setminus \{0\})).$
By Proposition 4.1, the set of singularities of $W(\mathscr{L}_{G}(\Sigma _*(G)))$ is the union of the critical value sets of $\mathbb{LH}^\pm_W$ which is the union of unfolded lightcone focal sets
$
\mathbb{LF}^+_W\cup \mathbb{LF}^-_W.
$
Therefore, we have shown the following proposition.
\begin{Pro} Let $(W,\mathcal{S})$ be a world hyper-sheet in $\R^{n+1}_1$  and $G:W\times (\R^{n+1}_1\setminus W)\lon \R$ the Lorentz distance squared function.
Then we have the graph-like legendrian unfolding
$\mathscr{L}_{G}(\Sigma _*(G))\subset  J^1(\R^{n+1}_1,I)$ such that $$W(\mathscr{L}_{G}(\Sigma _*(G)))=\mathbb{LH}^+_W(U\times I\times (\R\setminus \{0\}))\cup \mathbb{LH}^-_W(U\times I\times (\R\setminus \{0\})).$$
\end{Pro}
We write 
\[
\mathbb{LH}^\pm_{(W,\mathcal{S})}=\mathbb{LH}^\pm_W(U\times I\times (\R\setminus \{0\})).
\]
We also call $\mathbb{LH}^+_{(W,\mathcal{S})}\cup \mathbb{LH}^-_{(W,\mathcal{S})}$ an {\it unfolded light sheet} of $(W,\mathcal{S}).$
On the other hand, we have the corresponding Lagrangian submanifold $\Pi (\mathscr{L}_{G}(\Sigma _*(G)))\subset T^*\R^{n+1}_1.$
We now consider the natural question what are the caustic $C_{\mathscr{L}_{G}(\Sigma _*(G))}$ and  the Maxwell set $M_{\mathscr{L}_{G}(\Sigma _*(G))}$?
Moreover, are there any meanings of $C_{\mathscr{L}_{G}(\Sigma _*(G))}$ and $M_{\mathscr{L}_{G}(\Sigma _*(G))}$ in Physics?
%%%%%%%%%%%%%%%%%%%%%%
\section{Caustics of world hyper-sheets}
%%%%%%%%%%%%%%%%%%%%%%
\par
 In \cite{Bousso, Bousso-Randall} Bousso and Randall gave an idea of caustics of world sheets in order to define the notion of holographic domains.
The family of light sheets $\{\mathbb{LH}^\pm_{\mathcal{S}_{t}}(U\times\{t\})\times\R\}_{t\in J}$ sweeps out a region in $\R^{n+1}_1.$
A {\it caustic} of a world sheet is the union of the sets of critical values of light sheets along momentary spaces $\{\mathcal{S}_t\}_{t\in I}.$
A {\it holographic domain} of the world sheet is the region where the light-sheets sweep out until {\it caustics}.
So this means that the boundary of the holographic domain consists the caustic of the world sheet.
The set of critical values of the light sheet of a momentary space is the lightlike focal set of the momentary space.
Therefore the notion of caustics in the sense of Bousso-Randall is formulated as follows:
 {\it Caustics of a world sheet} $(W,\mathcal{S})$ 
 are defined to be
\[
\displaystyle{C^\pm(W,\mathcal{S})=\bigcup _{t\in I} \mathbb{LF}^\pm_{\mathcal{S}_{t}}}=\pi_1(\mathbb{LF}^\pm _{(W,\mathcal{S})}),
\]
where $\pi_1:\R^{n+1}_1\times I\lon \R^{n+1}_1$ is the canonical projection.
We call $C^\pm(W,\mathcal{S})$ {\it BR-caustics} of $(W,\mathcal{S}).$
We write that $C(W,\mathcal{S})=\pi_1(\mathbb{LF}^+ _W\cup \mathbb{LF}^- _W)$ and call it a {\it total BR-caustic} of $(W,\mathcal{S}).$
By definition, we have $\Sigma (W(\mathscr{L}_{G}(\Sigma _*(G)))=\mathbb{LF}^+ _{(W,\mathcal{S})}\cup \mathbb{LF}^- _{(W,\mathcal{S})},$ so that we have the following proposition.
\begin{Pro} Let $(W,\mathcal{S})$ be a world hyper-sheet in $\R^{n+1}_1$ and $G:U\times I\times (\R^{n+1}_1\setminus W)\lon \R$ 
the Lorentz distance squared function. Then we have
$C(W,\mathcal{S})=C_{\mathscr{L}_{G}(\Sigma _*(G))}.$
\end{Pro}
\par
In \cite{Bousso,Bousso-Randall} the authors did not consider the Maxwell set of a world sheet.
However, the notion of Maxwell sets plays an important role in the cosmology which has been called a {\it crease set} by Penrose (cf. \cite{Siino,Penrose}).
Actually, the topological shape of the event horizon is determined by the crease set of light sheets.
Here, we write $M(W,S)=M_{\mathscr{L}_G(\Sigma _*(G))}$ and call it a {\it BR-Maxwell set} of the world sheet $(W,\mathcal{S}).$ 
\par
Let $\bX_i:(U_i\times I_i, (\ou_i,t_i))\lon (\R^{n+1}_1,p_i)$, $(i=1,2)$ be germs of timelike embeddings such that
$(W_i,\mathcal{S}_{i})$ are world hyper-sheet germs, where $W_i=\bX_i(U).$ 
For $\bm{\lambda}_i=\mathbb{LH}^+_{
\mathcal{S}_{i}}(p_i,\ou_i)$ or $\bm{\lambda}_i=\mathbb{LH}^-_{
\mathcal{S}_{i}}(p_i,\ou_i),$ let $G_i :(U_i\times I_i\times (\R^{n+1}_1\setminus W_i) ,(\ou_i,t_i,\bm{\lambda}_i))\lon \R$
be Lorentz distance squared function germs.
We also write that $g_{i,\bm{\lambda}_i}(\ou,t)=G_i(\ou,t,\bm{\lambda}_i).$
Since $$W(\mathscr{L}_{G_i}(\Sigma _*(G_i)))=\mathbb{LH}^+_{(W_i,\mathcal{S}_i)}\cup \mathbb{LH}^-_{(W_i,\mathcal{S}_i)},$$
we can apply Theorem 6.2 and Corollary 6.3 to our case. Then we have the following theorem.

\begin{Th} Suppose that the set of critical points of $\overline{\pi}|_{\mathscr{L}_{G_i}(\Sigma _*(G_i))}$ are nowhere dense for $i=1,2$, respectively.
Then the following conditions are equivalent{\rm :}
\par\noindent
{\rm (1)} $(\mathbb{LH}^+ _{(W_1,\mathcal{S}_1)}\cup \mathbb{LH}^- _{(W_1,\mathcal{S}_1)}, \bm{\lambda}_1)$ and $(\mathbb{LH}^+ _{(W_2,\mathcal{S}_2)}\cup \mathbb{LH}^- _{(W_2,\mathcal{S}_2)}, \bm{\lambda}_2)$ are $S.P^+$-diffeomorphic,
\par\noindent
{\rm (2)} $\mathscr{L}_{G_1}(\Sigma _*(G_1))$ and $\mathscr{L}_{G_2}(\Sigma _*(G_2))$ are $S.P^+$-Legendrian equivalent,
\par\noindent
{\rm (3)} $\Pi(\mathscr{L}_{G_1}(\Sigma _*(G_1)))$ and $\Pi(\mathscr{L}_{G_2}(\Sigma _*(G_2))$ are Lagrangian equivalent.
\end{Th}
\par
We remark that conditions (2) and (3) are equivalent without any assumptions (cf. Theorem 6.2).
Moreover, if we assume that $\mathscr{L}_{G_i}(\Sigma _*(G_i))$ are $S.P^+$-Legendrian stable, then we can apply Proposition 5.3 and Theorem 6.6 and show the following theorem.
\begin{Th} Suppose that $\mathscr{L}_{G_i}(\Sigma _*(G_i))$ are $S.P^+$-Legendrian stable for $i=1,2,$ respectively.
Then the following conditions are equivalent{\rm :}
\par\noindent
{\rm (1)} $(\mathbb{LH}^+ _{(W_1,\mathcal{S}_1)}\cup \mathbb{LH}^- _{(W_1,\mathcal{S}_1)}, \bm{\lambda}_1)$ and $(\mathbb{LH}^+ _{(W_2,\mathcal{S}_2)}\cup \mathbb{LH}^- _{(W_2,\mathcal{S}_2)}, \bm{\lambda}_2)$ are $S.P^+$-diffeomorphic,
\par\noindent
{\rm (2)} $\mathscr{L}_{G_1}(\Sigma _*(G_1))$ and $\mathscr{L}_{G_2}(\Sigma _*(G_2))$ are $S.P^+$-Legendrian equivalent,
\par\noindent
{\rm (3)} $\Pi(\mathscr{L}_{G_1}(\Sigma _*(G_1)))$ and $\Pi(\mathscr{L}_{G_2}(\Sigma _*(G_2))$ are Lagrangian equivalent,
\par\noindent
{\rm (4)} $g_{1,\bm{\lambda}_1}$ and $g_{2,\bm{\lambda}_2}$ are $S.P$-$\mathcal{K}$-equivalent,
\par\noindent
{\rm (5)} $SK(\overline{W}_1, TLC(\mathcal{S}_{t_0},\bm{\lambda}_1)\times I;(p_1,t_0))=
SK(\overline{W}_2, TLC(\mathcal{S}_{t_0},\bm{\lambda}_2)\times I;(p_2,t_0)).$ 
\end{Th}
\par
By definition and Proposition 8.1, we have the following proposition.
\begin{Pro} If $\Pi(\mathscr{L}_{G_1}(\Sigma _*(G_1)))$ and $\Pi(\mathscr{L}_{G_2}(\Sigma _*(G_2))$ are Lagrangian equivalent,
then total BR-caustics $C(W_1,\mathcal{S}_1)$, $C(W_2,\mathcal{S}_2)$ and BR-Maxwell sets $M(W_1,\mathcal{S}_1)$, $M(W_2,\mathcal{S}_2)$ are diffeomorphic as set germs, respectively.
\end{Pro}

\section{World sheets in $\R^3_1$}
In this section we consider world sheets in the $3$-dimensional Minkowski space-time as an example.
Let $(W, \mathcal{S})$ be a world sheet in $\R^3_1$, which is parameterized by
a timelike embedding $\bm{\Gamma} :J\times I\lon \R^{3}_1$ such that $\mathcal{S}_t=\bm{\Gamma} (J\times \{t\})$ for $t\in I.$
In this case we call $\mathcal{S}_t$ a {\it momentary curve}.
We assume that $s\in J$ is the arc-length parameter. 
Since $W$ is a timelike surface, we have the unit pseudo-normal vector field $\bm{n}(s,t)$ of $W$ in $\R^{3}_1.$
Then  
$\bm{t}(s,t)=\bm{\gamma} _t'(s)$ is the unit spacelike tangent vector of $\mathcal{S} _t$.
We also define $\bm{b}(s,t)=\bm{n} (s,t)\wedge \bm{t}(s,t)$, which is the unit timelike normal vector of $\mathcal{S} _t$ in $TW$.
We choose the orientation of $\mathcal{S} _t$ such that $\bm{b}(s,t)$ is {\it future directed} (i.e. $\langle \bm{e}_0,\bm{b}(s,t)\rangle <0$).
Therefore,  
$\{\bm{b}(s,t),\bm{n}(s,t),\bm{t}(s,t)\}$ is a {pseudo-orthonormal frame} along $W. $ 
On this moving frame, we can show the following
{\it Frenet-Serret type formulae} for $S_t$:
\[
	\setlength\arraycolsep{2pt}
	\left\{
	\begin{array}{ccl}
	\displaystyle{\frac{\partial \bm{b}}{\partial s}(s,t)} &=& \tau_g(s,t) \bm{n}(s,t)-\kappa_g(s,t) \bm{t}(s,t), \\
    \displaystyle{\frac{\partial \bm{n}}{\partial s}(s,t)} &=& \tau_g(s,t) \bm{b}(s,t)-\kappa_n(s,t) \bm{t}(s,t) , \\
	\displaystyle{\frac{\partial \bm{t}}{\partial s}(s,t)} &=& -\kappa_g(s,t) \bm{b}(s,t) + \kappa_n(s,t) \bm{n}(s,t),
	\end{array}
	\right. 
\]
where $\kappa_g(s,t) = \langle \frac{\partial \bm{t}}{\partial s}(s,t), \bm{b}(s,t) \rangle,$
	$\kappa_n(s,t) = \langle \frac{\partial \bm{t}}{\partial s}(s,t), \bm{n}(s,t) \rangle,$
	$\tau _g(s,t)= \langle \frac{\partial \bm{b}}{\partial s}(s,t),\bm{n}(s,t)\rangle.$
We call $\kappa _g(s,t)$ a {\it geodesic curvature}, $\kappa _n(s,t)$ a {\it normal curvature} and $\tau _g(s,t)$ a {\it geodesic torsion} of $\mathcal{S}_t$ respectively.
It is known that
\par
   1) $\mathcal{S}_{t_0}$ is an asymptotic curve of $W$ if and only if $\kappa _n(s,t_0)\equiv 0,$
   \par
   2) $\mathcal{S}_{t_0}$ is a geodesic of $W$ if and only if $\kappa _g(s,t_0)\equiv 0,$
   \par
   3) $\mathcal{S}_{t_0}$ is a line of principal curvature of $W$ if and only if $\tau _g(s,t_0)\equiv 0.$
   \par\noindent
Then $\bm{b}(s,t_0)\pm \bm{n}(s,t_0)$ are lightlike. It follows that we have the light sheets
$\mathbb{LS}^\pm_{\mathcal{S}_{t_0}}:J\times\{t_0\}\times \R \lon \R^3_1$ along $\mathcal{S}_{t_0}$ defined by
$
\mathbb{LS}^\pm_{\mathcal{S}_{t_0}}((s,t_0),u)=\bm{\Gamma} (s,t_0)+u(\bm{b}(s,t_0)\pm \bm{n}(s,t_0)).
$
Here, we use the notation $\mathbb{LS}^\pm_{\mathcal{S}_{t_0}}$ instead of $\mathbb{LH}^\pm_{\mathcal{S}_{t_0}}$
because the images of these mappings are lightlike surfaces.
We adopt $\bm{n}^T=\bm{b}$ and $\bm{n}^S=\bm{n}.$ By the Frenet-Serret type formulae, we have
\[
\frac{\partial (\bm{n}^T\pm \bm{n}^S)}{\partial s}(s,t)=\frac{\partial (\bm{b}\pm \bm{n})}{\partial s}(s,t)=\tau_g(s,t))(\bm{n}\pm\bm{b})(s,t)-(\kappa _g(s,t)\pm \kappa _n(s,t))\bm{t}(s,t).
\]
Therefore, we have $\kappa ^\pm (\mathcal{S}_t)(s,t)=\kappa _g(s,t)\pm \kappa _n(s,t).$
It follows that 
\[
\mathbb{LF}^\pm _{\mathcal{S}_{t_0}}=\left\{\bm{\Gamma}(s,t_0)+\frac{1}{\kappa _g(s,t_0)\pm \kappa _n(s,t_0)}(\bm{b}\pm \bm{t})(s,t_0)\bigm | s\in J, \kappa _g(s,t_0)\pm \kappa _n(s,t_0)\not= 0\right\}.
\]
We consider the Lorentz distance squared function $G:J\times I\times \R^3_1\lon \R$.
In \cite{Izu2014}, by the calculation of the first and second derivative of $G$ with respect to $s$, we have introduced an invariant defined by
\[
\sigma^\pm (s,t)=((\kappa _n\pm\kappa _g)\tau _g\mp(\kappa_n'\pm\kappa_g'))(s,t).
\]
A geometrical meaning of these invariants is given as follows \cite{Izu2014}.
\begin{Pro}
The following conditions are equivalent{\rm :}
\par\noindent
{\rm (1)} $\sigma^\pm (s,t_0)\equiv 0$,\par\noindent
{\rm (2)} $\{\bm{\lambda}^\pm\}=\mathbb{LF}^\pm _{\mathcal{S}_{t_0}}$ is a point,
\par\noindent
{\rm (3)} $\mathcal{S}_{t_0}\subset LC_{\bm{\lambda}^\pm}.$
\end{Pro}
Moreover, as an application of the matrix criterion for $\mathcal{R}$-versality in \cite[Section 6.10]{Bru-Gib}, we have shown the following proposition in \cite{Izu2014}.
\begin{Pro}
{\rm (1)} The light sheet $\mathbb{LS}^\pm_{\mathcal{S}_{t_0}}(J\times\{t_0\}\times \R)$ at $\bm{\lambda}_0\in \mathbb{LF}^\pm_{\mathcal{S}_{t_0}}$ is local diffeomorphic to the cuspialedge $\bm{CE}$ if $\sigma ^\pm (s_0,t_0)\not= 0,$
\par
{\rm (2)} The light sheet $\mathbb{LS}^\pm_{\mathcal{S}_{t_0}}(J\times\{t_0\}\times \R)$ at $\bm{\lambda }_0\in \mathbb{LF}^\pm_{\mathcal{S}_{t_0}}$ is local diffeomorphic to the swallowtail $\bm{SW}$ if $\sigma ^\pm (s_0,t_0)= 0$ and $\partial \sigma ^\pm/\partial s (s_0,t_0)\not= 0.$
\par\noindent
Here, $\bm{CE}=\{(u,v^2,v^3)\ |\ (u,v)\in \R^2\}$, $\bm{SW}=\{(3u^4+u^2v,4u^3+2uv,v)\ |\ (u,v)\in \R^2\}$.
\end{Pro} 
\par
On the other hand, we now classify $S.P^+$-Legendrian stable graph-like Legendrian unfoldings $\mathscr{L}_G(\Sigma _*(G))$ by
$S.P^+$-Legendrian equivalence. By Theorems 6.5 and 6.6, it is enough to classify $\overline{f}$ by $S.P$-$\mathcal{K}$-equivalence under the condition that
\[
\dim _{\R} \frac{{\cal E}_{1+1}}{\left\langle \frac{\partial \overline{f}}{\partial q},\overline{f} \right\rangle _{{\cal E}_{1+1}}+
\left\langle \frac{\partial \overline{f}}{\partial t} \right\rangle _{\R}} \leq 3.
\]
In \cite{Izudoc,Izu95} we have the following proposition.
\begin{Pro}
With the above condition, $\overline{f}:(\R\times\R,0)\lon (\R,0)$ with $\partial\overline{f}/\partial t(0)\not= 0$ is $S.P$-$\mathcal{K}$-equivalent to one of the following germs{\rm :}
\par
{\rm (1)} $q,$
\par
{\rm (2)} $\pm t\pm q^2,$
\par
{\rm (3)} $\pm t +  q^3,$
\par
{\rm (4)} $\pm t \pm q^4,$
\par
{\rm (5)} $\pm t + q^5.$
\end{Pro}

The infinitesimally $S.P^+$-$\mathcal{K}$-versal unfolding
$\mathcal{F}:(\R\times (\R^3\times \R),0)\lon (\R,0)$ of each germ in the above list is given as follows (cf. \cite[Theorem 4.2]{Izu95}):
\par
(1) $q$
\par
(2) $\pm t\pm q^2,$
\par
(3) $\pm t + q^3+x_0q,$
\par
(4) $\pm t\pm q^4+x_0q+x_1q^2,$
\par
(5) $\pm t+q^5+x_0q+x_1q^2+x_2q^3.$
\par\noindent

By Theorem 6.6, we have the following classification.
\begin{Th} Let $(W,\mathcal{S})$ be a world sheet in $\R^3_1$ parametrized by a timelike embedding $\bm{\Gamma} :J\times I\lon \R^{3}_1$
and $G:J\times I\times \R^3_1\lon \R$ be the Lorentz distance squared function of $(W,\mathcal{S}).$
Suppose that the corresponding graph-like Legendrian unfolding $\mathscr{L}_G(\Sigma _*(G))\subset J^1(\R^3_1,I)$ is $S.P^+$-Legendrian stable.
Then the germ of the unfolded light sheet $\mathbb{LH}^+_{(W,\mathcal{S})}\cup \mathbb{LH}^-_{(W,\mathcal{S})}$ at any point is $S.P^+$-diffeomorphic to one of the following set germs in 
$(\R^3\times \R,0)${\rm :}
\par
{\rm (1)} $\{(u,v,w),0)\ |\ (u,v,w)\in (\R^3,0)\ \},$

\par
{\rm (2)} $\{(-u^2,v,w),\pm 2u^3)\ |\ (u,v,w)\in (\R^3,0)\ \},$
\par
{\rm (3)} $\{(\mp 4u^3-2vu,v,w),3u^3\pm vu^2)\ |\ (u,v,w)\in (\R^3,0)\ \},$
\par
{\rm (4)} $\{((5u^4+2vu+3wu^2,v,w),\pm(4u^4+vu^2+2wu^3))\ |\ (u,v,w)\in (\R^3,0)\ \}.$
\end{Th}
\demo
For any $(s_0,t_0,\bm{\lambda}_0)\in J\times I\times \R^3_1,$ the germ of $\mathscr{L}_G(\Sigma _*(G))\subset J^1(\R^3_1,I)$ at $\bm{z}_0=\mathscr{L}_G(s_0,t_0,\bm{\lambda}_0)$
is $S.P^+$-Legendrian stable. It follows that the germ of $g_{\bm{\lambda}_0}$ at $(s_0,t_0)$ is $S.P$-$\mathcal{K}$-equivalent to one of the germs in the list of Proposition 9.3.
By Theorem 6.6, the graph-like Legendrian unfolding $\mathscr{L}_G(\Sigma _*(G))$ is $S.P^+$-Legendrian equivalent to the graph-like Legendrian unfolding 
$\mathscr{L}_{\mathcal{F}}(\Sigma_*(\mathcal{F}))$ where $\mathcal{F}$ is the infinitesimally $S.P$-$\mathcal{K}$-versal unfolding of one of 
the germs in the list of Proposition 9.3.
It is also equivalent to the condition that the germ of the graph-like wave front $W(\mathscr{L}_{\mathcal{F}}(\Sigma_*(\mathcal{F})))$ is $S.P^+$-diffeomorphic to the corresponding graph-like wave front of one of the normal forms.
For each normal form, we can obtain the graph-like wave front.
We only show that case (5) in Proposition 9.3.
In this case we consider $\mathcal{F}(q,x_0,x_1,x_2,t)=\pm t+q^5+x_0q+x_1q^2+x_2q^3.$
Then we have
\[
\frac{\partial \mathcal{F}}{\partial q}=5q^4+x_0+2x_1q+3x_1q^2,
\]
so that the condition $\mathcal{F}=\partial \mathcal{F}/\partial q=0$ is equivalent to the condition that
\[
x_0=-(5q^4+x_0+2x_1q+3x_1q^2),\ t_0=\pm(4q^5+x_1q^2+2x_2q^3).
\]
If we put $u=q, v=x_0,w=x_1,$ then we have
\[
W(\mathscr{L}_{\mathcal{F}}(\Sigma_*(\mathcal{F})))=\{((-(5u^4+2vu+3wu^2),v,w),\pm(4u^4+vu^2+2wu^3))| (u,v,w)\in (\R^3,0)\}.
\]
It is $S.P^+$-diffeomorphic to the set germ of (4).
We have similar calculations for other cases.
We only remark here that we obtain the germ of (1) for both the germs of (1) and (2) in Proposition 9.3.
Since $W(\mathscr{L}_{\mathcal{F}}(\Sigma_*(\mathcal{F})))=\mathbb{LH}^+_{(W,\mathcal{S})}\cup \mathbb{LH}^-_{(W,\mathcal{S})},$ this completes the proof.
\enD
As a corollary, we have a local classification of BR-caustics in this case.
\begin{Co} With the same assumption for the world sheet $(W,\mathcal{S})$ as Theorem 9.4, the total BR-caustic $C(W,\mathcal{S})$ of $(W,\mathcal{S})$ at a singular point is locally diffeomorphic to the cuspidaledge $\bm{CE}$ or the swallowtail $\bm{SW}$.
\end{Co}
\demo
The total BR-caustic $C(W,\mathcal{S})$ of $(W,\mathcal{S})$ is the set of the critical values of $\pi_1\circ \overline{\pi}|_{\mathscr{L}_G(\Sigma _*(G))}.$
Therefore, it is enough to calculate the set of critical values of $\pi_1\circ \overline{\pi}|_{\mathscr{L}_{\mathcal{F}}(\Sigma_*(\mathcal{F}))}$
for each normal form $\mathcal{F}$ in Proposition 9.3.
For the germ (5) in Proposition 9.3, by the proof of Theorem 9.4 we have
$$
\Sigma_*(\mathcal{F})=\{(u,5u^4+2vu+3wu^2,v,w)\in (\R\times (\R^3\times \R),0)|(u,v,w)\in (\R^3,0)\}.
$$
It follows that
$$
\pi_1\circ \overline{\pi}\circ \mathscr{L}_{\mathcal{F}}(u,5u^4+2vu+3wu^2,v,w)=(5u^4+2vu+3wu^2,v,w).
$$
Then the Jacobi matrix of $f(u,v,w)=(5u^4+2vu+3wu^2,v,w)$ is
\[
J_f=\begin{pmatrix}
20u^3+2v+6wu & 0 & 0 \cr
2u & 1 & 0 \cr
3u^2 & 0 & 1 \cr
\end{pmatrix},
\]
so that the set of critical values of $f$ is  given by
\[
\{(-(15u^4+3wu^2),-10u^3-3wu,w)\in (\R^3,0)\ |\ (u,w)\in (\R^2,0)\}.
\] 
If we consider a linear transformation $\psi:(\R^3,0)\lon \R^3,0)$ defined by
$\psi (x_0,x_1,x_2)=(-\frac{1}{5}x_0,-\frac{2}{5}x_1,\frac{3}{5}x_2),$ then
we have
$\psi(-(15u^4+3wu^2),-10u^3-3uw,w)=(3u^4+\frac{3}{5}wu^2,4u^3+\frac{6}{5}wu,\frac{3}{5}w).$
If we put $U=u, V=\frac{3}{5}w,$ then we have $(3U^4+VU^2,4U^3+2VU,V),$ which is the parametrization of $\bm{SW}.$
By the arguments similar to the above, we can show that the set of critical values of $\pi_1\circ \overline{\pi}|_{\mathscr{L}_{\mathcal{F}}(\Sigma_*(\mathcal{F}))}$ is a regular surface for
the germ (3) and is diffeomorphic to $\bm{CE}$ for the germ (4) in Proposition 9.3, respectively.
This completes the proof.
\enD
\begin{Rem}{\rm
Since a world sheet $(W,\mathcal{S})$ is a timelike surface in $\R^3_1,$ we can define the evolute of $(W,\mathcal{S})$ by
\[
Ev_{(W,\mathcal{S})}=\displaystyle{\bigcup _{i=1}^2\left\{\bm{X}(u,t)+\frac{\displaystyle{1}}{\displaystyle{\kappa _i(u,t)}}\bm{n}^S(u,t)\ |\ u\in U,\kappa _i(u,t)\not= 0\ \right\}},
\]
where $\kappa _i(s,t)$ $(i=1,2)$ are the principal curvatures of $W$ at $p=\bm{X}(u,t).$
The evolute of a timelike surface has singularities in general.
The singularities of the evolute of a generic timelike surface in $\R^3_1$ is classified into $\bm{CE}$, $\bm{SW},$ $\bm{PY}$ or $\bm{PU}$ (cf., \cite{Farid}),
where $\bm{PY}=\{(u^2-v^2+2uv,-2uv+2uw,w)|w^2=u^2+v^2\}$ is the {\it pyramid} and $\bm{PU}=\{(3u^2+wv,3v^2+wu,w)|w^2=36uv\}$ is the {\it purse}.
So the singularities of BR-caustics of world sheets are different from those of the evolutes of timelike surfaces.
Since the singularities of BR-caustics are only corank one singularities, the pyramid and the purse never appeared in general.
Moreover, the normal line of a timelike surface is a spacelike line, so that it is not a ray in the sense of the relativity theory.
Therefore, the evolute of a timelike surface in the Minkowski space-time is not a caustic in the sense of Physics. 
}
\end{Rem}

\begin{flushleft}
\textsc{
Shyuichi Izumiya
\\ Department of Mathematics
\\ Hokkaido University
\\ Sapporo 060-0810, Japan} 
\par
e-mail: izumiya@math.sci.hokudai.ac.jp
\end{flushleft}
\end{document}